\magnification=1200\overfullrule=0pt\tolerance=2000
\documentstyle{amsppt}
\hcorrection{.25in}
\define\({\left(}
\define\){\right)}
\define\[{\left\[}
\define\]{\right]}
\redefine\epsilon{\varepsilon}
\define\R{\Bbb R^n}
\define\s{\Bbb S^{n-1}}

\define\la{\lambda}

\redefine\phi{\varphi}
\define\eb{\Cal B}

\topmatter \title A note on the $M^*$-limiting convolution body \endtitle
\author antonis tsolomitis\endauthor
\address Department of Mathematics,
        231 W.18th Avenue,
	Columbus, OH 43210 \endaddress
\curraddr 
Mathematical Sciences Research Institute, 1000 Centennial
Road, Berkeley, CA 94720
\endcurraddr 
\email atsol\@math.mps.ohio-state.edu\endemail

\abstract We introduce the mixed convolution bodies of two convex
symmetric bodies. We prove that if the boundary of a body $K$ is smooth enough
then as $\delta$ tends to $1$ the $\delta$--$M^*$--convolution body
of $K$ with itself tends to a multiple of the Euclidean ball after
proper normalization. On the other hand we show that the 
$\delta$--$M^*$--convolution body of the $n$--dimensional cube 
is homothetic to the unit ball of~$\ell_1^n$.
\endabstract

\endtopmatter
\document
\head 1. Introduction \endhead
\footnotetext""{Research at MSRI partially supported by NSF grant
DMS-9022140.}%

Throughout this note $K$ and $L$ denote convex symmetric bodies in $\R$.
Our notation will be the standard notation that can be found, for example,
in \cite{2} and \cite{4}.
 For $1\leq m\leq n$, $V_m (K)$ denotes the $m$--th mixed volume of
$K$ (i.e. mixing $m$ copies of $K$ with $n-m$ copies of the Euclidean
ball $\eb_n$ of radius one in $\R$).  Thus if $m=n$ then $V_n (K)=\text{vol}_n
(K)$
and if $m=1$ then $V_1 (K)=w(K)$ 
the mean width of $K$.

 For $0<\delta<1$ we define the $m$--th mixed $\delta$--convolution body of the convex symmetric bodies
$K$ and $L$ in $\R$:

\definition{Definition \rom{1.2} } The $m$--th mixed $\delta$--convolution body of
$K$ and $L$ is defined to be the set,
$$C_m (\delta;K,L)=
\{x\in\R\ :\ V_m  \( K\cap (x+L)\)\geq \delta V_m (K\cap L )\} .$$
\enddefinition
It is a consequence of Brunn--Minkowski inequality for mixed volumes
that these bodies are convex.

If we write $h(u)$ for the support function of $K$ in the direction 
$u\in\s$ then we have,
$$w(K)=2M_K^* =2 \int_{\s} h(u) d\nu (u),\tag1.1 $$
where $\nu$ is the Lebesgue measure of $\R$ restricted on $\s$ and 
normalized so that $\nu (\s )=1$.
In this note we study the limiting behavior of $C_1 (\delta;K,K)$
(which we will abbreviate with $C_1 (\delta)$)
as $\delta$ tends to $1$ and $K$ has a $C_+^2$ boundary. For simplicity
we will call $C_1 (\delta)$ the 
\lq\lq $\delta$--$M^*$--convolution body of $K$".

 We are looking for suitable $\alpha\in\Bbb R$ so that
the limit 
$$\lim_{\delta\rightarrow 1^- } \frac{C_1 (\delta)}{(1-\delta)^\alpha}
$$
exists (convergence in the Hausdorff distance). In this case we call
the limiting body \lq\lq the limiting $M^*$--convolution body of $K$".

We prove that for a convex symmetric body $K$ in $\R$
with $C_+^2$ boundary the limiting $M^*$--convolution body of $K$ 
is homothetic to the Euclidean ball. We also get a sharp estimate (sharp with
respect to the dimension $n$) of
the rate of the convergence of the $\delta$--$M^*$--convolution
body of $K$ to its limit. By $C_+^2$ we mean that the boundary of
$K$ is $C^2$ and that the principal curvatures of $bd(K)$ at every
point are all positive.


We also show that some smoothness condition on the boundary of $K$ is
necessary for this result to be true, by proving that the limiting
$M^*$--convolution body of the $n$--dimensional cube is homothetic to
the unit ball of $\ell_1^n$.

We want to thank Professor V.D.Milman for his encouragement and his guidance
to this research and for suggesting the study of mixed convolution bodies.

\head 2. The case \lq\lq $bd(K)$ is a $C_+^2$ manifold"\endhead

In this section we prove the following:
\proclaim{Theorem \rom{2.1}} Let $K$ be a convex symmetric body in $\R$ so that
$bd(K)$ is a $C_+^2$ manifold.
Then for all $x\in\s$ we have,
$$\bigg\vert \Vert x\Vert_{\frac{C_1 (\delta)}{1-\delta}} -\frac{c_n}{M_K^*}
\bigg\vert\leq C\frac{c_n}{M_K^*} 
\( M_K^* n (1-\delta) \)^2  ,\tag$2.1.1$ $$
where $c_n =\int_{\s} |\langle x,u\rangle | d\nu (u)
\sim 1/\sqrt{n}$ and $C$ is a constant 
independent of the dimension $n$. In particular,
$$\lim_{\delta\rightarrow 1^-} \frac{C_1 (\delta)}{1-\delta} =
\frac{M_K^*}{c_n}  \eb_{n} .$$
Moreover the estimate $(2.1.1)$ is sharp with respect to the dimension $n$
\endproclaim

By \lq\lq sharp" with respect to the dimension $n$ we mean that there are
examples (for instance the $n$--dimensional Euclidean ball) for which
the inequality $(2.1.1)$ holds true if \lq\lq $\leq$" is substituted
with \lq\lq $\geq$" and the constant $C$ is adjusted by a (universal) constant 
factor.

Before we proceed with the proof we will need to collect some standard
notation which can be found in \cite{4}.
We write $p:\ bd(K)\rightarrow \s$ for the Gauss map $p(x)=N(x)$ where
$N(x)$ denotes the unit normal vector of $bd(K)$ at $x$. $W_x$ denotes
the Weingarten map, that is, the differential of $p$ at the point $x\in bd(K)$.
$W_u^{-1}$ is the reverse Weingarten map at $u\in\s$ and the eigenvalues 
of $W_x$ and 
$W_u^{-1}$ are respectively the principal curvatures and principal 
radii of curvature of the manifold $bd(K)$ at $x\in bd(K)$ and $u\in\s$. 
We write $\Vert W\Vert$ and $\Vert W^{-1}\Vert$ for the quantities:
$\sup_{x\in bd(K)} \Vert W_x\Vert$ and $\sup_{u\in\s} \Vert W_u^{-1} \Vert$
respectively. These quantities are finite since the manifold $bd(K)$ is 
assumed to be $C_+^2$.

 For $\la\in\Bbb R$ and $x\in\s$ we write $K_\la$ for the set 
$K\cap (\la x +K )$. 
$p_\la^{-1} :\ \s\rightarrow bd(K_\la )$ is the reverse Gauss map, that is,
the affine hyperplane $p_\la^{-1} (u)+[u]^\perp$ is tangent to $K_\la$
at $p_\la^{-1} (u)$. The normal cone of $K_\la$ at $x$ is denoted
by $N(K_\la ,x)$ and similarly for $K$. The normal cone is a convex set (see
\cite{4}). Finally $h_\la$ will denote the support function of $K_\la$.

\demo{Proof}
Without loss of generality we may assume that both the $bd(K)$ and $\s$ are
equipped with an atlas whose charts are functions which are Lipschitz,
their inverses are Lipschitz and they all have the same Lipschitz
constant $c>0$. 

Let $x\in\s$ and $\la=\frac1{\Vert x\Vert_{C_1 (\delta)}}$; hence
$\la x\in bd\( C_1 (\delta)\)$ and 
$$M_{K_\la }^* =\delta M_K^*  .\tag$2.1.2$ $$
We estimate now $M_{K_\la}^*$. Let $u\in\s$. We need to compare $h_\la (u)$
and $h(u)$. Set $Y_\la =bd(K)\cap bd(\la x+K)$. 

{\it Case 1}.  $p_\la^{-1} (u)\notin Y_\la$.

In this case it is easy to see that
$$ h_\la (u)= h(u)-|\langle \la x,u\rangle|.$$

{\it Case 2}. $p_\la^{-1} (u)\in Y_\la$.

Let $y_\la =p_\la^{-1} (u)$ and $y'_\la =y_\la -\la x\in bd(K)$.
The set $N(K_\la , y_\la )\cap\s$ defines a curve $\gamma$ which we assume 
to be parametrized on $[0,1]$ with $\gamma (0)=N(K, y_\la)$ and $\gamma (1)=
N(K, y'_\la )$. We use the inverse of the Gauss map $p$ to 
map the curve $\gamma$ to a curve $\tilde{\gamma}$ on $bd(K)$ by setting
$\tilde{\gamma}=p^{-1} \gamma$. The end points of $\tilde{\gamma}$
are $y_\la$ (label it with $A$) and $y'_\la$ (label it with $B$).
Since $u\in\gamma$ we conclude that the point $p^{-1} (u)$ belongs
to the curve $\tilde{\gamma}$ (label this point by $\Gamma$).
Thus we get:
$$0\leq h(u)-h_\la (u)=|\langle \vec{A\Gamma}, u\rangle |.$$
It is not difficult to see that the cosine of the angle of the vectors
$\vec{A\Gamma}$ and $u$ is less than the largest principal curvature
of $bd(K)$ at $\Gamma$ times $|\vec{A\Gamma} |$, the length of the vector $\vec{A\Gamma}$.
Consequently we can write,
$$0\leq h(u)-h_\la (u)\leq \Vert W\Vert |\vec{A\Gamma}|^2 .$$
In addition we have,
$$\align
|\vec{A\Gamma}|&\leq \text{length}\( \tilde{\gamma}|_A^\Gamma \)\leq
\text{length} \(\tilde{\gamma}|_A^B \)\\
&=\int_0^1 |d_t \tilde{\gamma} |dt=\int_0^1 |d_t p^{-1} \gamma |dt\\
&\leq\Vert W^{-1} \Vert \text{length} (\gamma)\leq\frac2{\pi} \Vert
W^{-1} \Vert |p(y_\la )-p(y'_\la )|,
\endalign
$$
where $|\centerdot |$ is the standard Euclidean norm.
Without loss of generality we can assume that the points $y_\la$ and
$y'_\la$ belong to the same chart at $y_\la$. Let $\phi$ be the chart
mapping $\Bbb R^{n-1}$ to a neighborhood of $y_\la$ on $bd(K)$ and 
$\psi$ the chart mapping $\Bbb R^{n-1}$ on a neighborhood of 
$N(K, y_\lambda )$ in $\s$. We assume, as we may, that
the graph of $\gamma$ is contained in the range of the chart $\psi$.
It is now clear from the above series of inequalities that
$$|\vec{A\Gamma}|\leq c_0 \Vert W^{-1} \Vert |\psi^{-1} p \phi (t)-
\psi^{-1} p \phi (s)|,$$
where $t$ and $s$ are points in $\Bbb R^{n-1}$ such that $\phi (t)=y_\la$
and $\phi (s)=y'_\la$ and $c_0 >0$ is a universal constant.
Now the mean value theorem for curves gives,
$$\align
|\vec{A\Gamma} |&\leq C\Vert W^{-1} \Vert \Vert W\Vert |t-s|\\
&\leq  C\Vert W^{-1} \Vert \Vert W\Vert |y_\la -y'_\la |\\
&=C\Vert W^{-1} \Vert \Vert W\Vert \la ,
\endalign
$$
where $C$ may denote a different constant every time it appears.
Thus we have,
$$0\leq h(u)-h_\la (u)\leq C\Vert W\Vert\(\Vert W^{-1} \Vert \Vert W\Vert\)^2
\la^2 .$$

Consequently,
$$
 \int_{\s\setminus p_\la (Y_\la )} \( h(u)-|\langle \la x,
u\rangle|\) d\nu (u) + \int_{p_\la (Y_\la )} \( h(u)-C \la^2 \) d\nu (u)$$
$$\leq M_{K_\la}^* =\delta M_K^* \leq$$
$$\int_{\s\setminus p_\la (Y_\la )} \( h(u)-|\langle \la x,
u\rangle|\) d\nu (u) + \int_{p_\la (Y_\la )}  h(u)d\nu (u) ,$$
where $C$ now depends on $\Vert W\Vert$ and $\Vert W^{-1} \Vert$.

Rearranging and using $c_n$ for the quantity $\int_{\s} |\langle x,u\rangle |
d\nu (u)$ and the fact $\la=1/\Vert x\Vert_{C_1 (\delta)}$ we get:
$$\bigg\vert \Vert x\Vert_{\frac{C_1 (\delta)}{1-\delta}} -\frac{c_n}{M_K^*}
\bigg\vert\leq \frac{c_n}{M_K^*} \(\frac{\int_{p_\la (Y_\la )}
|\langle x,u\rangle|d\nu (u)}{c_n} + C\la\frac{\mu\(p_\la (Y_\la )\)}{c_n} \).$$

We observe now that for $u\in p_\la (Y_\la ), |\langle x,u\rangle|\leq
\text{length} (\gamma)/2 \leq \Vert W\Vert \la$. Using this in the last 
inequality and the fact that $p_\la (Y_\la )$ is a \lq\lq band"
around an equator
of $\s$ of width at most $\text{length} (\gamma )/2$ we get,
$$\align
\bigg\vert \Vert x\Vert_{\frac{C_1 (\delta)}{1-\delta}} -\frac{c_n}{M_K^*}
\bigg\vert&\leq \frac{c_n}{M_K^*} C n\la^2 \\
&\leq \frac{c_n}{M_K^*} C n\frac{(1-\delta)^2 }{\Vert x\Vert_{
\frac{C_1 (\delta)}{1-\delta}}^2 }.
\endalign
$$

Our final task is to get rid of the norm that appears on the right side
of the latter inequality. Set $$T=\frac{\Vert x\Vert_{C_1 (\delta) /1-\delta }}{c_n /M_K^* } .$$
We have shown that $$T^2 \vert T-1 \vert \leq C
\frac{M_K^*}{c_n}  n (1-\delta )^2 .$$
If $T\geq 1$ then we can just drop the factor $T^2$ and we are done.
If $T<1$ we write $T^2 |T-1|$ as $\(1-(1-T)\)^2 (1-T)$ and we consider 
the function $$f(x)=(1-x)^2 x\ :\ (-\infty ,\frac13 )\rightarrow\Bbb R.$$
This function is strictly increasing thus invertible on its range, that is,
$f^{-1}$ is well defined and increasing in $(-\infty ,\frac4{27} )$.
Consequently if $$C\frac{M_K^*}{c_n} n(1-\delta)^2 \leq \frac4{27} ,
\tag$2.1.3$ $$
we conclude that,
$$\align
0\leq 1-T &\leq f^{-1} \( C\frac{M_K^*}{c_n} n(1-\delta)^2 \)\\
&\leq C\frac{M_K^*}{c_n} n(1-\delta)^2 .
\endalign
$$
The last inequality is true since the derivative of $f^{-1}$ at zero is $1$.
Observe also that the convergence is \lq\lq essentially realized"
after $(2.1.3)$ is satisfied.\qed

\enddemo

We now proceed to show that some smoothness conditions on the boundary of
$K$ are necessary, by proving that the limiting 
$M^*$--convolution body of the $n$--dimensional cube is homothetic to the
unit ball of $\ell_1^n$. In fact we show that the 
$\delta$--$M^*$--convolution body of the cube is already homothetic to
the unit ball of $\ell_1^n$.

\proclaim{Example \rom{2.3} }
Let $P=[-1,1]^n$. Then for $0<\delta<1$ we have,
 $$C_1 (P)=\frac{C_1 (\delta;P,P)}{1-\delta}
=n^{3/2} \text{vol}_{n-1} (\s) B_{\ell_1^n} .$$
\endproclaim

\demo{Proof}
Let $x=\sum_{j=1}^n x_j e_j$ where $x_j \geq 0$ for all $j=1,2,\dots ,n$
and $e_j$ is the standard basis of $\R$. Let $\la>0$ be such that
$\la x\in bd\( C_1 (\delta)\)$. Then,
$$P\cap (\la x+P)=\{ y\in\R :\ y=\sum_{j=1}^n y_i e_i , -1+\la x_j \leq
y_j \leq 1 , j=1,2,\dots ,n \}.$$
The vertices of $P_\la =P\cap (\la x+P)$ are the points $\sum_{j=1}^n
 \alpha_j e_j$ where $\alpha_j$ is either $1$ or $-1+\la x$ for all $j$.
Without loss of generality we can assume that $-1+\la x_j <0$ for all
the indices $j$. Put $\hbox{sign}\,\alpha_j =\alpha_j /|\alpha_j |$ when
$\alpha_j \neq 0$ and $\hbox{sign}\,0=0$. Fix a sequence of $\alpha_j$'s so that
the point $v=\sum_{j=1}^n \alpha_j e_j$ is a vertex of $P_\la$.
Clearly, $$N\( P_\la , v\)=N\(P, \sum_{j=1}^n (\hbox{sign}\,\alpha_j)e_j \).$$
If $u\in\s\cap N(P_\la ,v)$ then,
$$h_\la (u)=h (u)-\bigg\vert \langle \sum\limits_{j=1}^n (\alpha_j -
\hbox{sign}\,\alpha_j )e_j ,u\rangle \bigg\vert .$$
If $\hbox{sign}\,\alpha_j =1$ then $\alpha_j -\hbox{sign}\,\alpha_j =0$ otherwise $\alpha_j
-\hbox{sign}\,\alpha_j =\la x .$

Let $\Cal A \subseteq \{1,2,\dots ,n \}$. Consider the 
\lq\lq$\Cal A$--orthant":
$$\Cal O_{\Cal A} =\{y\in\R :\ \langle y, e_j \rangle <0 , \ \text{if}\ 
j\in\Cal A\ \text{and}\ \langle y,e_j \rangle \geq 0\ \text{if} \ j\notin\Cal A
\} .$$
Then $\Cal O_{\Cal A} =N\( P, \sum_{j=1}^n (\hbox{sign}\,\alpha_j )e_j \)$ if and only
if $\hbox{sign}\,\alpha_j =1$ exactly for every $j\notin\Cal A$.
Hence,
$$h_\la (u)=h(u)-\bigg\vert \langle\sum\limits_{j\in\Cal A} \la x_j e_j ,u
\rangle \bigg\vert,$$ 
for all $u\in\Cal O_{\Cal A} \cap \s$.
Hence using the facts $M_{P_\la}^* =\delta M_P^*$ and $\la=1/\Vert x\Vert_{
C_1 (\delta ) }$ we get,
$$\Vert x\Vert_{\frac{C_1 (\delta)}{1-\delta}} =-\frac1{M_P^*}
\sum\limits_{\Cal A\subseteq \{ 1,2,\dots , n\} } \sum\limits_{j\in\Cal A}
x_j \int_{\Cal O_{\Cal A} \cap \s } \langle e_j ,u\rangle d\nu (u) ,$$
which gives the result since 
$$\int_{\Cal O_{\Cal A} \cap \s } \langle e_j ,u\rangle d\nu (u)=
\frac1{2^{n-1}} \int_{\s} |\langle e_1 ,u\rangle| d\nu (u).$$\qed

\enddemo

\head{References}\endhead
\widestnumber\no{1}
\ref
\no 1
\by Kiener, K\.
\paper Extremalit\"at von Ellipsoiden und die Faltungsungleichung von Sobolev
\jour Arch\. Math.
\vol 46
\yr 1986 \pages 162--168\endref

\ref
\no 2
\by Milman, V\. and Schechtmann, G\.
\book Asymptotic theory of finite dimensional normed spaces
\publ Springer Lecture Notes
\vol 1200 \yr 1986
\endref

\ref
\no 3
\by Schmuckenschl\"ager, M\.
\paper The distribution function of the convolution square of a convex symmetric
body in $\Bbb R ^n$
\jour Israel Journal of Mathematics
\vol 78 \yr 1992 \pages 309--334
\endref

\ref
\no 4
\by Schneider, R\.
\book Convex bodies: The Brunn--Minkowski theory
\publ Encyclopedia of Mathematics and its Applications, Cambridge University Press
\vol 44
\yr 1993
\endref

\ref
\no 5
\by Tsolomitis, A\.
\paper Convolution bodies and their limiting behavior
\jour Duke J. Math. (submitted)
\endref

\enddocument